\documentclass[sanskrit]{article}
in \oddsidemargin 0.25true in \evensidemargin 0.25true in
\usepackage{epsfig}

\usepackage{enumerate}

\usepackage{graphicx}
\usepackage{CJKutf8}
\usepackage{devanagari}

\usepackage{amsmath}

\usepackage{lipsum}

\begin{document}

\newcommand{\arcsech}{\mbox{arcsech}}
\newcommand{\bc}{\begin{center}}
\newcommand{\ec}{\end{center}}
\newcommand{\be}{\begin{equation}}
\newcommand{\ee}{\end{equation}}
\newcommand{\bea}{\begin{eqnarray}}
\newcommand{\eea}{\end{eqnarray}}
\newcommand{\bean}{\begin{eqnarray*}}
\newcommand{\eean}{\end{eqnarray*}}
\newcommand{\bd}{\begin{description}}
\newcommand{\ed}{\end{description}}

\newcommand{\caln}{\cal N}
\newcommand{\caloi}{{\cal OI}}
\newcommand{\hatv}{{\bf v}}
\newcommand{\hatrho}{\hat{{\bf \rho}}}
\newcommand{\hatlam}{\lambda}
\newcommand{\tillam}{\tilde{\lambda}}
\newcommand{\tT}{\tilde{T}}
\newcommand{\hT}{\hat{T}}
\newcommand{\bT}{{\bf T}}
\newcommand{\bA}{{\bf A}}
\newcommand{\ep}{\epsilon}
\newcommand{\cd}{{\cal D}}
\newcommand{\hp}{\hat{T}}
\newcommand{\epf}{\epsilon_f}
\newcommand{\Nn}{N_{net}}
\newcommand{\Ns}{N_{sph}}
\newcommand{\Nd}{N_{dat}}
\newcommand{\lng}{\langle}
\newcommand{\rng}{\rangle}
\newcommand{\vep}{\epsilon}
\newcommand{\vph}{\varphi}
\newcommand{\vp}{\varphi}
\newcommand{\bfk}{{\bf k}}
\newcommand{\bfr}{{\bf r}}
\newcommand{\bfN}{{\bf N}}
\newcommand{\bfa}{{\bf a}}
\newcommand{\bfb}{{\bf b}}
\newcommand{\bfc}{{\bf c}}
\newcommand{\bfH}{{\bf H}}
\newcommand{\bfK}{{\bf K}}
\newcommand{\bfL}{{\bf L}}
\newcommand{\bfM}{{\bf M}}
\newcommand{\bfx}{{\bf x}}
\newcommand{\bC}{{\bf C}}
\newcommand{\bp}{{\bf{p}}}
\newcommand{\bfe}{{\bf e}}
\newcommand{\hbt}{\hat{\overline{T}}}
\newcommand{\uvr}{{\bf \hat{r}}}
\newcommand{\uvn}{{\bf \hat{n}}}
\newcommand{\bw}{\overline{w}}
\newcommand{\bx}{\mbox{\boldmath $\xi$}}
\newcommand{\bnu}{\mbox{\boldmath $\nu$}}
\newcommand{\my}{M_{yr}}

\newcommand{\bvt}{\begin{verbatim}}
\newcommand{\evt}{\end{verbatim}}
\newcommand{\dbar}{\vrule width 0.25truein height 1.0pt depth -0.6pt}

\renewcommand{\baselinestretch}{2.0}

{\raggedleft
{\bf\large
OFFER (One-Figure-Facilitates-Every-Relationship)
for Trigonometric Functions }

}
\vskip 0.5cm

\renewcommand{\baselinestretch}{1.2}

{\raggedleft 
{\large SAMUEL S.P. SHEN }  \\
San Diego State University, San Diego, CA 92182\\
{\tt sshen@sdsu.edu}\\
ORCID: 0000-0002-4535-4684

\vskip 0.0cm

}
\vskip 0.5cm

\thispagestyle{empty}
\renewcommand{\baselinestretch}{1.2}

\noindent   {\bf Abstract:} This article describes an approach that uses a single figure to illustrate all six trigonometric functions, their fundamental identities, and the law of cosines. We name this method OFFER (one-figure-facilitates-every-relationship) for the convenience of mathematics teaching. The article also describes the SPORT (story-picture-observe-review-tell) approach to effective mathematics learning. 


\section{Introduction}
The expression $\sin x$ is a trigonometric function, also known as a trig function for short. Does $\sin x$ mean the {\bf sin}  committed by someone named $x$? Of course, not. Nonetheless, the strange names and notations of trigonometric functions
have frustrated numerous students from high school to college, due to infinitely many trig identities and their variations. While instructors remember many identity formulas and can solve problems with ease, students do not have confidence in memorizing them because they believe these formulas are unrelated to their daily lives and do not form simple pictures in their minds. Even worse, the only picture some students have of trig is a nightmare! Grandparents do not even know the meaning of $\sin x$, and parents have forgotten what they learned, while elder siblings are helpless but empathetically share their own nightmares and curses. Many students feel that trig functions are simply tricky, make them unworthy, torture them, and treacherously block them from entering science or engineering. 

\begin{figure}[h]
\centering
\includegraphics[width=4in]{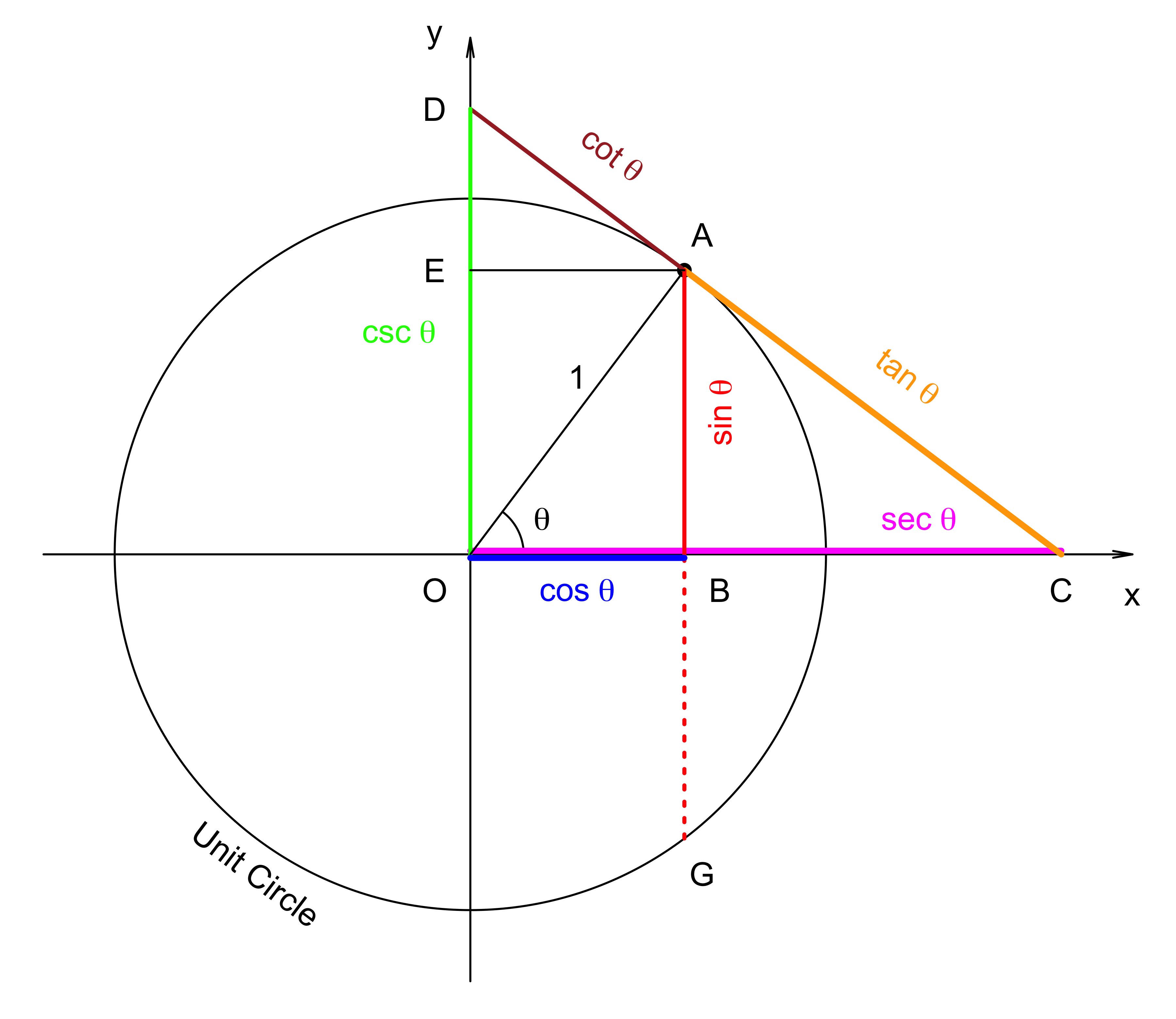}
\caption{All six trigonometric functions illustrated by a single figure. The R and Python code for generating this figure is available for free download [7]. Also see the code availability section of this paper.}
\label{fig01}
\end{figure}

With the hope of reducing students' frustration and increasing their confidence, we provide materials to help teach or learn trig functions effectively, efficiently, and with depth and insight. We propose to use the story-picture-observe-review-tell (SPORT) pedagogy. Our approach is to use a single figure to explain all the definitions, fundamental identities, and their proofs of trig functions.  This approach can be summarized as ``one-figure-facilitates-every-relationship (OFFER)." See Fig. \ref{fig01} that gives the geometric definitions of all six trigonometric functions. We hope to use this pictorial approach to dispel students' nightmares and to make them enjoy learning mathematics. 

SPORT may be regarded as a pedagogy of mathematics education, derived from learning-by-doing. SPORT assumes that students should first learn materials useful to them, substantiated by stories. Almost all the useful mathematics definitions, theorems, formulas, and algorithms have fantastic {\bf S}tories behind them. Each story must have a {\bf P}ictorial illustration. Hand-drawing such a picture can enhance the student's understanding of the mathematical theory, for example, drawing Fig. \ref{fig01} by hand. Teachers can inspire students to {\bf O}bserve the picture(s) and lead them to reinvent the theory in classrooms. Students can independently redraw the figures, observe them, and reach their own conclusions and extensions beyond textbooks. This {\bf S-P-O} process can be 
{\bf R}eviewed or repeated many times to enhance learning. Ancient Chinese educator Confucius said, "To learn the new, review the old."  Finally, students can master the materials and {\bf T}ell stories to others. This completes the entire SPORT learning cycle. 

The SPORT mathematics education, in nature, requires instructors to first focus on the most fundamental and most useful materials. one-figure-facilitates-every-relationship (OFFER) is a concrete approach to using the SPORT pedagogy. As an example of OFFER, this article uses only one figure to explain and derive all the fundamental trigonometric formulas. In this way, we can liberate students from massive textbooks and lead them to see the insight, power, depth, and beauty of mathematics. 


The purpose of this article is not to provide you with a textbook to teach trigonometry, but rather to recommend that you consider using geometry, intuition, figures, and common sense to teach or learn. This paper may also spark wide-ranging discussions among school teachers and university professors and, consequently, help modernize the teaching of trigonometry in the era of the digital and AI (artificial intelligence) economy. Therefore, this paper does not strictly follow the usual style of technical writing; rather, it uses the style of a conversation with you. 

We provide convenient R and Python code that allows both instructors and teachers to quickly reproduce and modify the single figure that explains all trigonometric functions and their fundamental identities. Although none of our materials is entirely new, the concise and fairly complete summary of OFFER for trigonometry with computing support might be useful to many.  

The next section will be the story of India's invention of sine and cosine functions. Section 3 contains the OFFER materials: use a single figure to explain the definitions of trigonometric functions and their fundamental identities, classified into two categories: identities based on right triangles and those based on oblique triangles. Section 4 describes the derivative of the sine function again using a single figure. Section 5 discusses the SPORT learning-by-doing pedagogy for mathematics education. Summary and discussion are in Section 6. Computer code availability is in Section 7.

\section{Story: Trig functions were invented in India around 500 AD}

When I learned trig functions in a middle school in a remote Chinese village, I was amazed at how my teacher could remember so many formulas. I thought I would never be able to do it, even though I had the best memory in my class! He was known as the best mathematics teacher in the entire county, but he told us no stories about these mysterious functions and failed to teach us which of the infinitely many identities were fundamental and most useful. Nevertheless, he did show the magic of a unit circle, and I remembered that. However, I could not remember most formulas, nor could I derive them, despite being a top student at the school taught by a top teacher. 

In fact, trigonometric functions could be taught in a much simpler way, with much deeper insight, by telling the story of how they were invented and spread around the world, and by drawing the corresponding figure. 

\begin{figure}[h]
\centering
\includegraphics[width=4in]{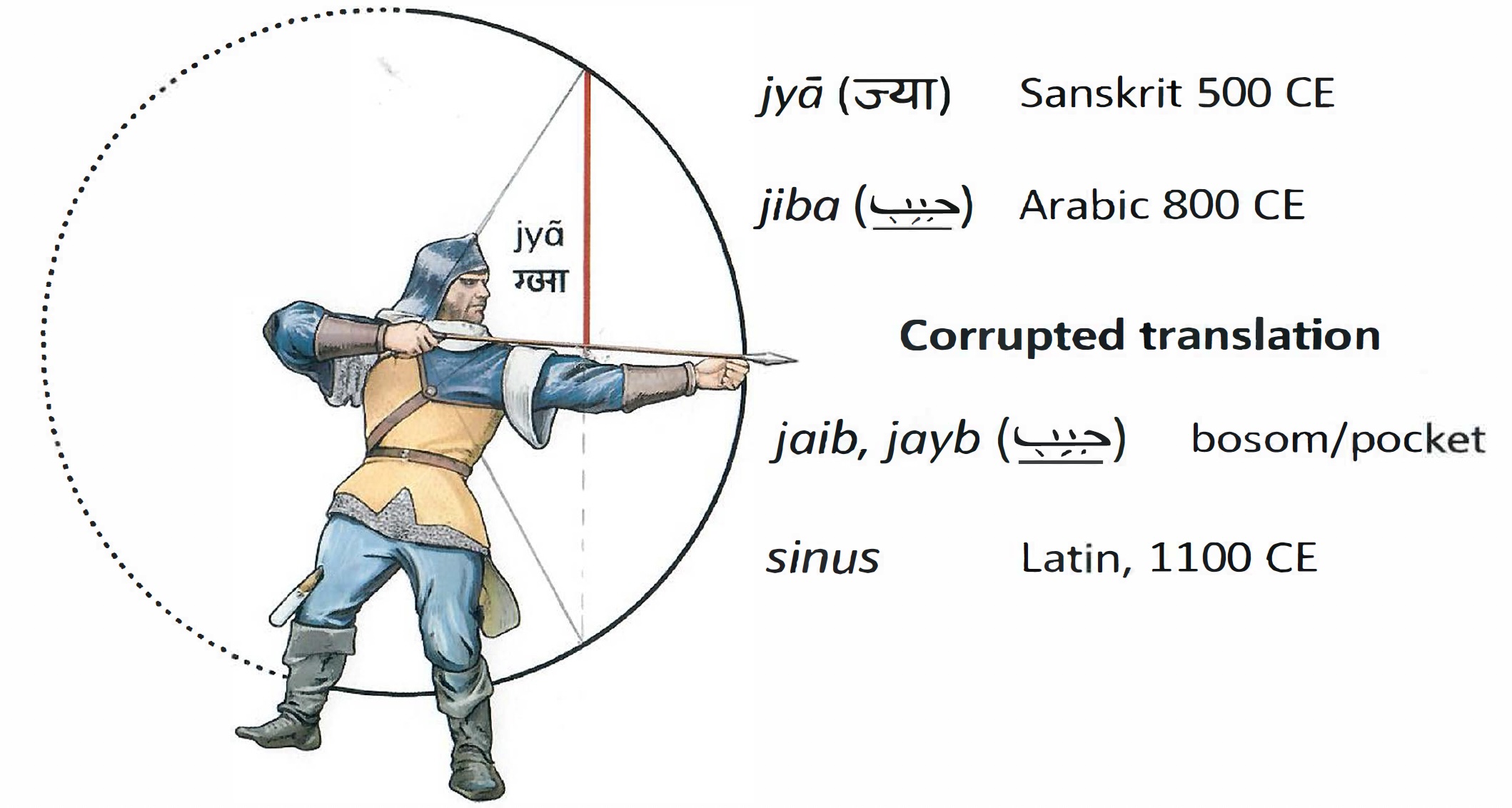}
\caption{Indian mathematician Aryabhata used half-chord, denoted by {\tt jya} in the figure, to define today's sine function. Credit: This figure was produced by ChatGPT following a series of prompts from the author.}
\label{fig02}
\end{figure}

Trigonometric functions, particularly the sine and cosine,  were invented by the Indian mathematician and astronomer Aryabhata (476 - 550 CE) [1]. He abstracted the bow-arrow shape in archery into a triangular geometry as shown in Fig. \ref{fig02}. 

Aryabhata regarded the rigid bow of an archery set as a part of a circle, and the elastic string as a chord. When the string is pulled, the upper half-string and the arrow form an angle $\theta$, as shown in Fig. \ref{fig02}. The original position of the upper half-string is called the half-chord, {\tt ardha-jya}, the Latin spelling of the Sanskrit term {\dn ardhajyA}. Here, {\tt ardha} means half, and {\tt jya} means chord. In Sanskrit mathematics writing, {\tt jya} was used, and {\tt ardha} was omitted. Aryabhata included this in his 499 CE book entitled "{\it Aryabhatiya}", written in Sanskrit [1, 2]. In the book, he computed a table of 
{\tt jya} values for different angles, the first known sine table [1]. For example, he had $\sin(\pi/48) = 0.0654, \sin(11\pi/48) = 0.6594$, which are very close to computer-generated values 0.0654 and 0.6593.

The meaning of {\tt jya}, i.e., the half-chord, got distorted in a series of translations from Sanskrit to Arabic, Latin, and English. The distorted English translation is the one used today in our books, classrooms, and computer programs. 

Sanskrit {\tt jya} was phonetically transcribed to Arabic as {\tt jiba} around 800 CE. 
In written Arabic, short vowels are often omitted. In this case, {\tt jiba} became {\tt j-b}. 
Later Arabic readers and the 1100 CE European translators mistakenly regarded the vowelless {\tt jb} as 
the commonly used word {\tt jaib}, since {\tt jiba} was a rare Arabic word, perhaps because it was a technical term.
This misreading led to a very different meaning, {\tt jaib} meaning ``bosom" of a female body, ``bay" in the bay window, ``pocket," ``bend," ``gulf," or ``fold," whose Latin correspondence is {\tt sinus}, meaning a curved surface. Later,  {\tt sinus} was anglicized into sine.  In the 1700s, Swiss mathematician Leonhard Euler popularized the function notation $\sin x$.  

Therefore, the later Arabic understanding and the European translation lost the original meaning of half-chord ({\tt ardha-jya}), or chord ({\tt jya}). However, in their 1859 Chinese translation book entitled ``Algebra and Calculus Step-by-Step" \begin{CJK*}{UTF8}{gbsn}(代微积拾级)\end{CJK*}, the Chinese mathematician Li Shanlan (last name Li) and the British Christian missionary Alexander Wylie respected the original meaning of ``jya."  They coined the Chinese name of ``sine function"  as {\tt zhengxian hanshu}  \begin{CJK*}{UTF8}{gbsn}(正弦函数)\end{CJK*}, literally meaning ``right chord function." Li and Wylie's right chord refers to the half-chord ``directly opposite" \begin{CJK*}{UTF8}{gbsn}(正对着)\end{CJK*} the angle. Although this meaning respects Aryabhata's original definition of the sine function, it would be a more accurate Chinese translation from Aryabhata's original meaning if Li and Wylie had used \begin{CJK*}{UTF8}{gbsn}半弦函数\end{CJK*}, meaning ``half-chord function."





\section{Trigonometric functions' definitions and identities based on a single figure}

\subsection{A figure with both right and oblique triangles}
\begin{figure}[h]
\centering
\includegraphics[width=4.5in]{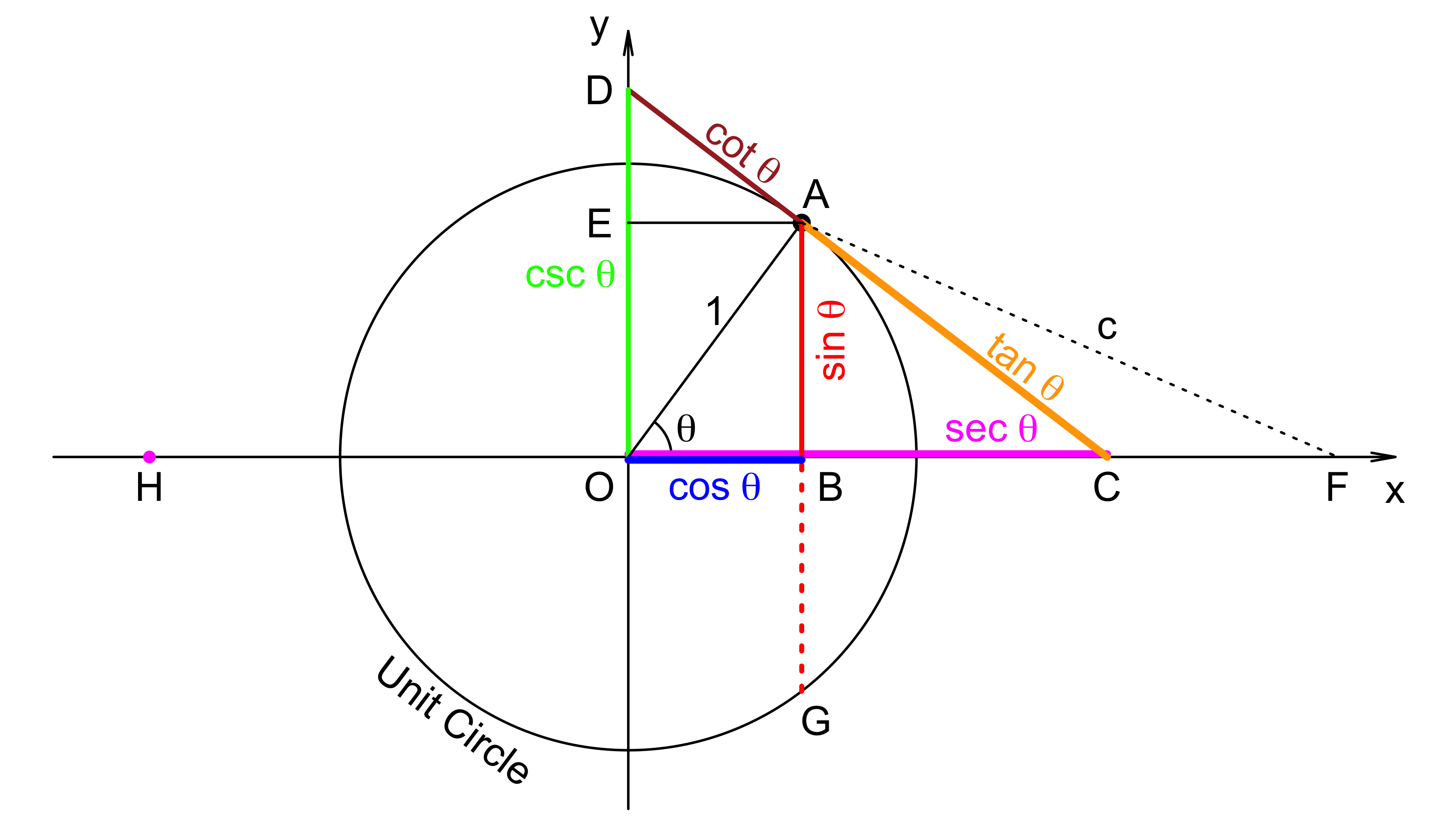}
\caption{A figure for the definition of trigonometric functions using the right triangles and for the proof of the law of cosines using an oblique triangle. }
\label{fig03}
\end{figure}

Figure \ref{fig03} is almost the same as Fig. \ref{fig01}, except that Fig. \ref{fig03} has an extra line AF that forms two oblique triangles, while all the triangles in Fig. \ref{fig01} are right. Showing something similar serves the purpose of review. Reviewing the old helps learn the new, in the line of Confucius' teaching. In publications, it is usually considered unprofessional to present the same figures or text more than once in an article. However, for learning, this is a good way. When students see something they have seen before, it may trigger their memory and prompt them to check for differences and consider why the author or instructor would do that. This is the psychology of curiosity. We human beings all have the capability of curiosity! However, it is still a challenge to inspire students to think and to spark their curiosity about exploration. 

\subsection{Definitions of trigonometric functions}
\subsubsection{Sine and cosine functions}
For the right triangle $\Delta OBA$ inside the unit circle in Fig. \ref{fig03}, Aryabhata defined the length of the half-chord $\overline{AB}$ (i.e., 
a half of the entire chord is $\overline{AG}$)
as $\sin\theta = AG$ and produced what is known as the first table of $\sin$ values for different angles. 
The length of $\overline{OB}$, regarded as the complementary half chord, is defined as $\cos\theta = OB$, i.e., 
complementary sine of the angle $\theta$. 

When the angle $\theta$ changes, $\sin \theta$ and $\cos \theta$ vary and are called sine and cosine functions of the angle $\theta$. For some special angles, the values of $\sin \theta$ and $\cos \theta$ can be inferred from Fig. \ref{fig03} and the above definitions. 
For example, 
\bea
&& \sin 0 = 0, ~~~~ \cos 0 = 1,\\
&& \sin \frac{\pi}{2} = 1, ~~~~ \cos \frac{\pi}{2} = 0, \\
&& \sin \frac{\pi}{4} = \frac{\sqrt{2}}{2}, ~~~~ \cos \frac{\pi}{4} = \frac{\sqrt{2}}{2}.
\eea

You may check the following properties of the trigonometric functions by 
drawing a figure like Fig. \ref{fig03}, but with $- \theta$, i.e., move point $A$ to its mirror position with respect to the $x$-axis. 
Your figure will imply the following
\be
\sin (-\theta) = - \sin \theta, ~~~~  \cos (-\theta) =  \cos \theta.
\ee




\subsubsection{Tangent and cotangent functions}
According to the Oxford English Dictionary, ``tangent" originates from the Latin {\tt tangent}, meaning ``touching," derived from the verb {\tt tangere} (``to touch"). It was first used in English in the late 1500s to describe the tangent line
in mathematics. The orange and brown line segments are on the tangent line to the unit circle at point $A$, determined by angle $\theta$. The distance between $A$ and the intersection point $C$ with the $x$-axis, which is the arrow line of Fig. \ref{fig02}, is defined as $\tan \theta = AC$, also denoted by  $\text{tg} \theta$ in some textbooks. Thus, the tangent function is the length of a tangent line segment. The name and the geometric meaning agree, unlike the sine or sinus function, whose English name and geometric meaning do not agree. 

Most textbooks define the tangent function using the following ratio
\be
\tan \theta = \frac{\sin \theta}{\cos \theta},
\ee
interpreted as the ratio of the opposite side $AB$ to the adjacent side $OB$. 
Compared to the visual line-segment definition, this ratio definition is abstract, algebra-based, and memory-dependent, and is not preferred in education because a ratio may not form a mental image for students. Actually, this ratio can be derived from the geometric definition, the length of the tangent line segment $\tan \theta = AC$. The derivation can be made from two similar right triangles $\Delta OAC$ and $\Delta OBA$, since $AC \perp AO$ and $AB \perp OC$. Two similar triangles have equal ratios of corresponding side lengths. Here, we use the legs of the right triangles (i.e., the two sides adjacent to the right angle):
\be
\frac{CA}{OA} = \frac{AB}{OB} ~~\longrightarrow ~~\frac{\tan \theta}{1}=\frac{\sin \theta}{\cos \theta}.
\ee
Thus, we may regard the ratio definition of $\tan \theta$ as a consequence of the intuitive geometric definition. 

The complementary tangent, denoted $\cot\theta$ or $\text{ctg}\theta$, is the length of the tangent line segment from $overline{AD}$, where $ D$ is the intersection of the tangent line with the $y$-axis. The line segment $AD$ may be regarded as being complementary to $overline{AC}$. Again, from this common-sense definition, you can derive the ratio consequence that appears as a definition in most textbooks and is used in many classrooms:
\be
\cot\theta =  \frac{\cos  \theta}{\sin \theta}.
\ee

Consider the two similar right triangles $\Delta OAD$ and $\Delta OBA$. The ratios of corresponding legs are 
\be
\frac{DA}{OA} = \frac{OB}{AB} ~~\longrightarrow ~~\frac{\cot \theta}{1}=\frac{\cos \theta}{\sin \theta}.
\ee
Namely,
\be
\cot \theta=\frac{\cos \theta}{\sin \theta}.
\ee

\subsubsection{Secant and cosecant functions}
 The word ``Secant" is derived from Latin ``secans," meaning ``cut." A secant line in mathematics is a line that cuts through 
 a curve. For the secant function in trigonometry, the ``cut" is related to the line segment $\overline{HC}$ on the$x-$axis that cuts the unit circle into two parts (see Fig. \ref{fig03}). 
The secant function $\sec\theta = OC$ is defined as the length of the segment $\overline{OC}$ and is half of the length of 
 the horizontal secant line segment $\overline{HC}$. Similar to the half-chord definition of $\sin \theta$, 
$\sec\theta$ uses only the right half length of the secant line segment $\overline{HC}$ as shown in Fig. \ref{fig03}. 

From this geometric definition and picture in mind, you can also derive the algebraic definition:
\be
\sec \theta = \frac{1}{\cos \theta}.
\label{eq10}
\ee
This ratio is used to define the secant function in most textbooks and classrooms, similar to the textbook definition of the tangent function. Again, this definition is abstract, requires students' memory, and is not preferred for teaching and learning. Instead, a better teaching approach should treat this expression as a consequence of the geometric definition of the secant function. You may draw by hand a figure similar to Fig. \ref{fig03} a couple of times, enhance your understanding of $\sec\theta$, and tell a story about the function.

The derivation for Eq. (\ref{eq10}) can be made by considering the two similar right triangles $\Delta OAC$ and $\Delta OBA$. The ratios of the corresponding two hypotenuses and two legs are
\be
\frac{OC}{OA} = \frac{OA}{OB} ~~\longrightarrow ~~\frac{\sec \theta}{1}=\frac{1}{\cos \theta}.
\ee

The complementary secant function is called the cosecant function, defined as the length of the vertical line segment $\overline{OD}$ and denoted by $\csc\theta$. The ``cut" that corresponds to the cosecant function is the secant line on the $y-$ axis 
that cuts the unit circle into a left half and a right half. 

With this geometric picture and definition in mind, you can also derive the algebraic definition
\be
\csc \theta = \frac{1}{\sin \theta}
\ee
from two similar right triangles $\Delta OAD$ and $\Delta OBA$. The ratios of the corresponding two hypotenuses and two legs are
\be
\frac{OD}{OA} = \frac{OA}{AB} ~~\longrightarrow ~~\frac{\csc \theta}{1}=\frac{1}{\sin \theta}.
\ee

\subsubsection{Summary of the six trig functions and an example of confusion}

In summary, the geometric definition of the six trigonometric functions is defined as follows:
\bea
&& \sin \theta = AB  ~~\mbox{(Length of a half-chord} ~ \overline{AB} ~), \label{eqdef1} \\
&& \cos \theta = OB ~~\mbox{(Length of a complementary half-chord} ~ \overline{OB} ~),\\
&& \tan \theta = AC ~~\mbox{(Length of a tangent line segment} ~ \overline{AC} ~), \\
&& \cot \theta = AD ~~\mbox{(Length of a complementary tangent line segment} ~ \overline{AD} ~),\label{eqdef2}\\
&& \sec \theta = OC ~~\mbox{(Length of a half secant line segment} ~ \overline{OC} ~), \\
&&  \csc \theta = OD ~~\mbox{(Length of a complementary half secant line segment} ~ \overline{OD} ~). \label{eqdef3}
\eea
Students can draw a diagram like Fig. \ref{fig03}, form a picture of the six trigonometric functions in their minds, have something concrete to think about when they want to, and tell stories about these functions with friends and family members aided by gestures and drawing. 

On the other hand, the ratio definition of the trigonometric functions is not preferred and may cause confusion.
For example, many students, even mature scientists, are often puzzled by the ratio definition of $\sec \theta$ and $\csc \theta$:
\be
\sec \theta = \frac{1}{\cos \theta} , ~~~~ \csc \theta = \frac{1}{\sin \theta}.
\ee
Why is the secant defined by cosine, while the cosecant is defined by sine? 
 This is against intuition and does not make sense! 
I have been asked this question many times, some out of curiosity, and others seeking answers to the questions from their children or students.  Their confusion may come from the pronunciation: {\tt SEE} in $\sec$ linked to {\tt KO} in $\cos$, and {\tt KO} in $\csc$ linked to {\tt SAI} in $\sin$. 

All this unnecessary confusion exists among many because they did not know the intrinsic geometric definition of the six trigonometric functions shown in Eqs. (\ref{eqdef1}) - (\ref{eqdef3}). The ``complementary" function of the secant is the cosecant, and the 
ratio of $\csc \theta = \frac{1}{\sin \theta}$ is only a derived result, which is a consequence, not a definition. Therefore, with the intrinsic geometric definition, this confusion may go away. 

Although it is not wrong mathematically to use the ratio $\csc \theta = 1/\sin \theta$ as a definition for the cosecant function, 
pedagogically it is not a good one. It is not a preferred way to teach or to learn trigonometry. The geometric definitions of the trigonometric functions have a clear advantage, as they are based on common sense, figures, and mental images. With the geometric definitions, in the classroom and in homework assignments, instructors may repeatedly ask students to hand-draw figures like Fig. \ref{fig03} and to explore their drawings to obtain additional results, such as the trigonometric identities shown in the next subsection. 

\subsection{Fundamental trigonometric identities}

\subsubsection{Trigonometric identities based on the Pythagorean theorem for right triangles}

\begin{figure}[h]
\centering
\includegraphics[width=3.5in]{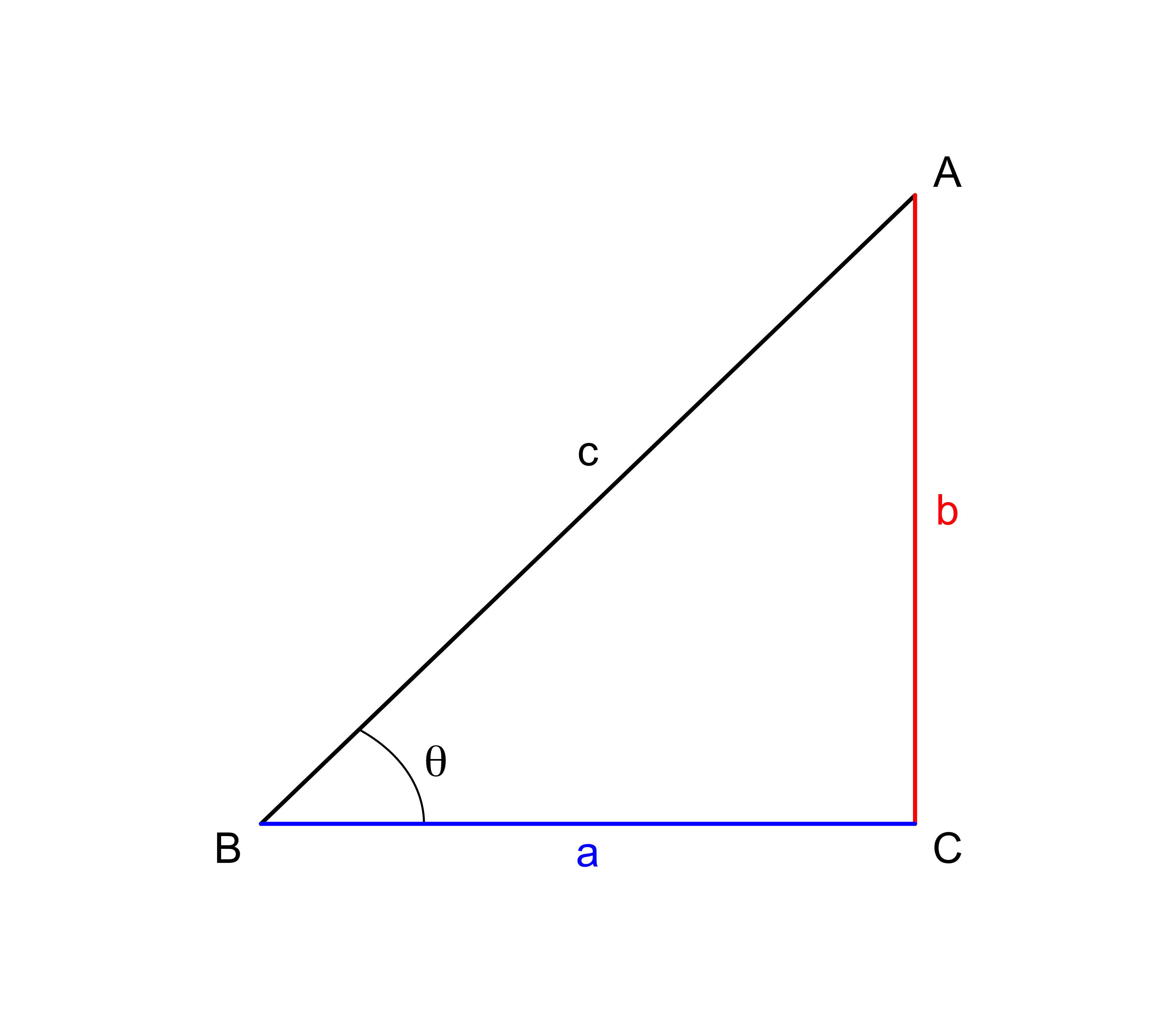}
\caption{A right triangle to support the expression (\ref{eq21}) of the Pythagorean theorem. }
\label{fig04}
\end{figure}

The Pythagorean theorem for a right triangle (Fig. \ref{fig04}) may be written as follows:
\be
a^2 + b^2 = c^2,
\label{eq21}
\ee
where $a$ and $b$ are two legs and $c$ is the hypotenus as shown in the figure. 

Each right triangle in Fig. \ref{fig03} corresponds to a Pythagorean theorem formula; in turn, it corresponds to a trigonometric identity. For example, the right triangle $\Delta OBA$ corresponds to $a = \cos\theta, b = \sin \theta$ and $c =1$ in Fig. \ref{fig04}, and Pythagorean theorem formula $a^2 + b^2 = c^2$ corrsponds to
\be
\sin^2 \theta + \cos^2 \theta  = 1^2 = 1.
\label{eq16}
\ee
This is the most popular and most fundamental trigonometric identity. 

To get an identity for $\tan \theta$, we can examine a right triangle with $\tan \theta$ as a side, say, 
$\Delta OAC$ in Fig. \ref{fig03}. This triangle corresponds to the following Pythagorean theorem
\be
1^2 + \tan^2 \theta = \sec^2 \theta.
\label{eq23}
\ee

This identity can also be derived from the fundamental identity (\ref{eq16}) by dividing both sides of the equation by
$\cos^2 \theta$:
\be
\frac{\sin^2 \theta + \cos^2 \theta}{\cos^2 \theta}  =  \frac{1}{\cos^2 \theta},
\label{eq24}
\ee
or
\be
\frac{\sin^2 \theta}{\cos^2 \theta} + 1  =  \frac{1}{\cos^2 \theta}.
\ee
This implies (\ref{eq23}).

This is an algebraic approach and is more abstract than the intuitive right-triangle approach for trigonometry beginners and those who lack confidence in their algebraic skills. Yet, most existing textbooks use the algebraic approach. We suggest minimizing the algebraic approach and using the triangle approach whenever possible, because it is easier to build a picture in mind with a triangle and to tell a story. 

Of course, the algebraic approach has its advantage of clean, direct, and efficient expressions. Students with a very good mathematics background may even prefer to start with the fast algebraic approach. In practical teaching, both approaches should be taught. However, regretably, most existing textbooks do not discuss the triangle approach at all. Therefore, few students can tell the geometric meaning of trig functions, let alone trig identities. This does not imply ignoring the algebraic approach. On the contrary, the important algebraic approach must be taught comprehensively, but only after the geometric triangle approach is taught. This helps students understand the geometric meaning of each algebraic step, form mental pictures, and tell stories. 

When students have formed good mental images, they can further explore the geometric meaning of the algebraic operations. For example, operation (\ref{eq24}) means 
\be
\frac{\mbox{AB}^2 + \mbox{OB}^2}{\mbox{OB}^2} = \frac{1}{\mbox{OB}^2}
\ee
The expression 
\bea
&&\frac{\mbox{OB}^2}{\mbox{OB}^2} =1^2 = \mbox{OA}^2,\\
&&\frac{\mbox{AB}^2}{\mbox{OB}^2} = \tan^2\theta = \mbox{CA}^2,\\
&&\frac{1}{\mbox{OB}^2} = \sec^2\theta = \mbox{CO}^2.
\eea
This implies that $\cos \theta$ is the dilation factor between the two similar triangles, $\Delta OBA$ and $\Delta OAC$. 
We may denote this relationship by 
\be
\Delta OBA \cong \Delta OAC \times \cos \theta.
\ee

Similar to the trig identity Eq. (\ref{eq23}) for $\tan \theta$, triangle $\Delta OAD$ leads to an identity for $\cot \theta$:
\be
1^2 + \cot^2 \theta = \csc^2 \theta.
\ee

The three examples above imply that any right triangle can yield a trigonometric identity. For example, triangle $\Delta ABC$ in Fig. \ref{fig03} leads to another identity:
\be
\sin^2\theta  + (\sec \theta - \cos \theta)^2 = \tan^2 \theta.
\ee
This complicated identity is not considered a fundamental identity and is rarely taught in classrooms. It is not surprising that it can also be algebraically derived from Eq. (\ref{eq23}).  

Would you like to identify a right triangle in Fig. \ref{fig03} and write its corresponding trigonometric identity? 

\subsubsection{Law of cosines for oblique triangles}
Now, let us consider an angle-and-side relationship for the oblique triangle $\Delta OAF$ in Fig. \ref{fig03}. We want to find out how side $\overline{AF}$ is related to two other sides $\overline{OA}$ and $\overline{OF}$, as well as the angle $\theta$ opposite to $\overline{AF}$. Let $a=OF$ be the length of the line segment $\overline{OF}$ and $c=AF$. The Pythagorean theorem for 
the right triangle $\Delta ABF$ yields
\be
AF^2  = AB^2  +  BF^2.
\label{eq20}
\ee
This can be written as
\be
c^2 = \sin^2\theta  + (a  - \cos \theta)^2 = \sin^2\theta + a^2 - 2a \cos \theta + \cos \theta^2 .
\ee
A simplification yields 
\be
c^2 = 1 + a^2 - 2a \cos \theta .
\label{eq24a}
\ee
This is called the law of cosines because the formula involves $\cos \theta$ in an oblique triangle. 
This law determines the opposite side of a given angle $\theta$ when the two other sides of the angle are also given.
The law of cosines applies to any triangle, oblique or right. For a right triangle, the law of cosines reduces to the Pythagorean theorem since $\cos(\pi/2) = 0$. Thus, the law of cosines may be regarded as a generalized Pythagorean theorem. 

If the unit-length side is stretched to length $r$, then multiply the previous equation by $r^2$ to obtain
\be
(rc)^2 = r^2 + (ra)^2 - 2r\times (ra) \cos \theta .
\ee
Let $s = ra$ and $t = rc$. Then, this equation becomes 
\be
t^2 = r^2 + s^2 - 2rs \cos \theta .
\ee
This is a formula for the law of cosines presented in most textbooks. 

The law of cosines implies that any triangle, not only right triangles, can lead to an identity involving trigonometric functions. 
For example, the law of cosines for the oblique triangle $\Delta ACF$ in Fig. \ref{fig03} implies the following identity:
\be
\mbox{AF}^2 = \tan^2\theta + \mbox{CF}^2 + 2  \times \mbox{CF} \times \tan\theta \sin \theta,
\label{eq38}
\ee
because $\cos(\pi/2 + \theta) = -\sin \theta$, which may be derived by drawing another figure like Fig. \ref{fig03}. If $\theta = \pi/4$ and $\mbox{CF} =1$, then formula (\ref{eq38}) can yield the 
following result
\be
\mbox{AF}^2 = 2 + \sqrt{2}, 
\ee
or 
\be
\mbox{AF} \approx 1.8478.
\ee

\subsection{Trigonometric functions for the addition and subtraction of two angles}

\subsubsection{The cosine difference formula}
\begin{figure}[h]
\centering
\includegraphics[width=3.5in]{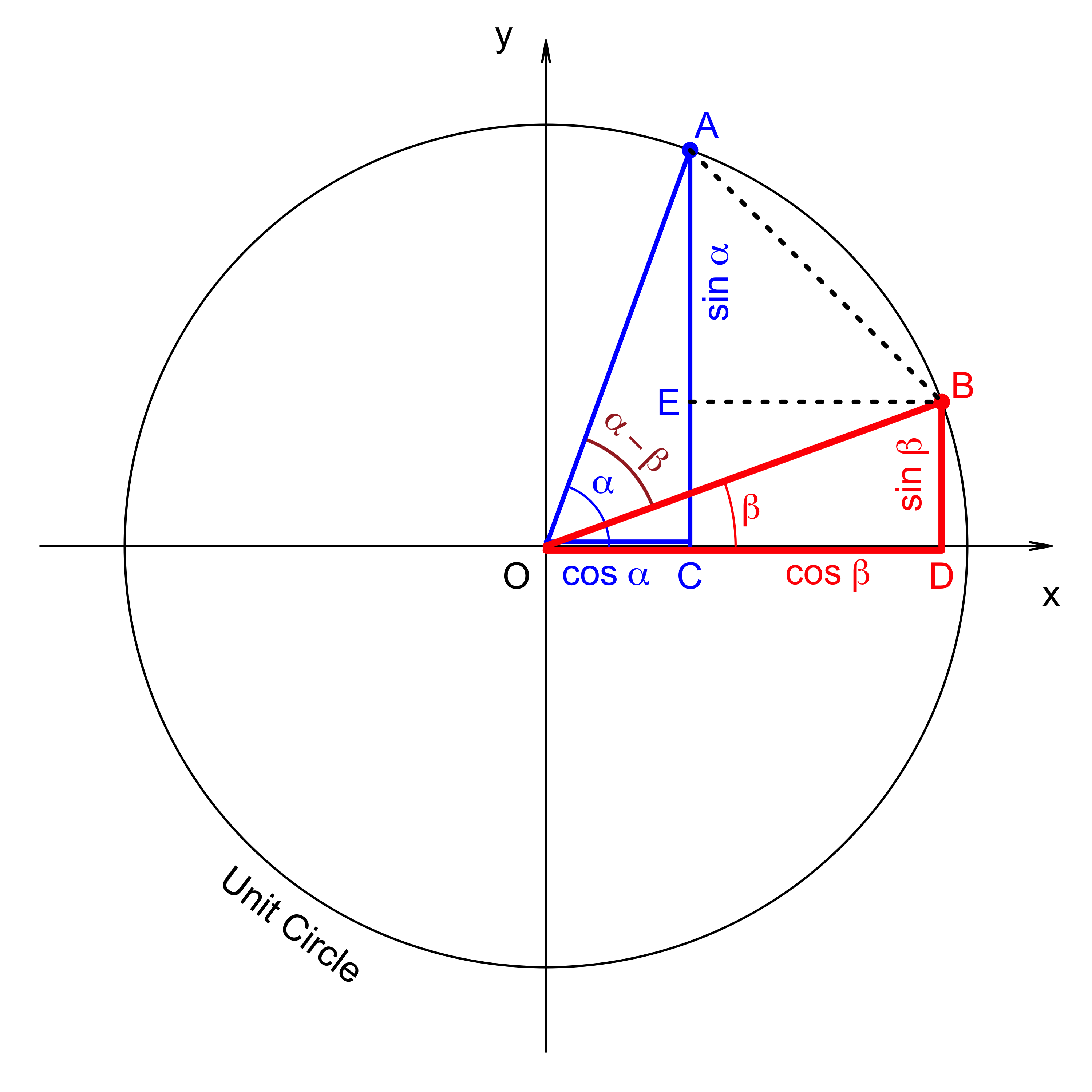}
\caption{A figure for the derivation of the cosine difference formula. }
\label{fig05}
\end{figure}

We explore the cosine of the difference of one angle minus another $\cos(\alpha - \beta)$, when
$\sin \alpha, \cos \alpha, \sin \beta$ and $\cos \beta$ are given.  
The result is called the cosine difference formula, or cosine subtraction formula.
 It has profound implications in both mathematics and science, including the mathematical 
 formulation for the concepts of stationarity, homogeneity, and isometry. This important formula can be derived using triangles in Fig.\ref{fig05}, the law of cosines, and the Pythagorean theorem. 

First, apply the law of cosines to the oblique triangle $\Delta AOB$ of Fig.\ref{fig05}:
\be
AB^2 = OA^2 + OB^2 - 2\times  OA\times
OB \times \cos(\alpha - \beta).
\label{eq25}
\ee
Here, $OA = OB = 1$ because $OA = OB$ is the radius of the unit circle. Equation (\ref{eq25}) yields
\be
AB^2 = 2 - 2\cos(\alpha - \beta).
\label{eq25a}
\ee

Next, apply the Pythagorean theorem to the right triangle $\Delta AEB$:
\be
AB^2 = BE^2 + AE^2.
\label{eq26}
\ee
Substitution of $BE = \cos\beta - \cos\alpha$ and $AE = \sin\alpha - \sin\beta$
into this formula yields
\bea
AB^2 &=& (\cos\beta - \cos\alpha)^2 + (\sin\alpha - \sin\beta)^2 \nonumber \\
&=&(\cos^2\beta - 2 \cos\beta \cos\alpha + \cos^2\alpha) + (\sin^2\beta - 2 \sin\beta \sin\alpha + \sin^2\alpha) \nonumber \\
&=&(\cos^2\alpha + \sin^2\alpha) + (\cos^2\beta + \sin^2\beta)- 2 \cos\beta \cos\alpha - 2 \sin\beta \sin\alpha \nonumber \\
&=&1 + 1- 2 \cos\beta \cos\alpha - 2 \sin\beta \sin\alpha \nonumber \\
&=& 2 - 2(\cos\alpha \cos\beta  + \sin\alpha \sin\beta ).
\label{eq27}
\eea

Comparison between Eqs. (\ref{eq25a}) and (\ref{eq27}) leads to 
\be
\cos(\alpha - \beta) = \cos\alpha \cos\beta  + \sin\alpha \sin\beta .
\label{eq28}
\ee
This is the cosine difference formula. Our derivation above is not much different from those in most textbooks. 
\footnote{The cosine difference formula has a popular 3D counterpart: the addition law for spherical harmonics. An infinite-dimensional counterpart can be demonstrated by the Fourier series for a homogeneous field. }

\subsubsection{The cosine addition formula}
The cosine addition formula can be easily derived when replacing $\beta$ by $-\beta$ in the cosine subtraction formula
\be
\cos(\alpha + \beta) = \cos(\alpha -(- \beta)) =  \cos\alpha \cos(-\beta)  + \sin\alpha \sin(-\beta).
\ee
This can be written as 
\be
\cos(\alpha + \beta) =  \cos\alpha \cos(\beta)  - \sin\alpha \sin(\beta),
\label{eq29}
\ee
since $\cos(-\beta) = \cos(\beta)$ and $\sin(-\beta) = -\sin(\beta)$.

\subsubsection{Addition and subtraction formulas for sine and tangent functions}

The sine addition formula can be derived in the following way:
\bea
&&\sin(\alpha + \beta) \nonumber \\
=&& \cos(\pi/2 - (\alpha + \beta))  \nonumber \\
=&& \cos((\pi/2 - \alpha) - \beta))  \nonumber \\
=&& \cos(\pi/2 - \alpha) \cos\beta  + \sin(\pi/2 - \alpha) \sin\beta  \nonumber \\
=&& \sin\alpha \cos\beta  + \cos\alpha \sin\beta.
\eea
So, the sine addition formula is
\be
\sin(\alpha + \beta) = \sin\alpha \cos\beta  + \cos\alpha \sin\beta.
\ee

The sine subtraction formula can be derived from the sine addition formula with $-\beta$
\be
\sin(\alpha -\beta) = \sin(\alpha + (-\beta) )= \sin\alpha \cos(-\beta)  + \cos\alpha \sin(-\beta).
\ee
This can be written as 
\be
\sin(\alpha -\beta) = \sin\alpha \cos\beta  - \cos\alpha \sin\beta.
\ee

The tangent addition formula can be derived by the ratio
\bea
\tan(\alpha -\beta) &&= \frac{\sin(\alpha -\beta)}{\sin(\alpha -\beta)} \nonumber \\
&&= \frac{\sin\alpha \cos\beta  - \cos\alpha \sin\beta}{\cos\alpha \cos\beta  + \sin\alpha \sin\beta } .
\eea
Division of the numerator and denominator by $\cos\alpha \cos\beta$ yields
\be
\tan(\alpha -\beta) =  \frac{\tan \alpha - \tan \beta}{1 + \tan \alpha  \tan \beta}. 
\ee

This formula is rarely used in applications, except for teaching.  For trigonometry beginners, I suggest that you ignore these kinds of materials that are Only-Useful-for-Teaching (OUT), at least the first couple of rounds of your study of a subject. In traditional precalculus and calculus textbooks and classrooms, many topics, such as the function $y = x^x$ and its derivative, are OUT materials. It is a more productive and effective learning experience if you first focus on Relevant, Useful, and Modern (RUM) materials. Only after you have formed clear mental images and feel comfortable with relevant algebraic operations, you can explore the more complex trig identities, including the OUT materials. By then, you can further ask what geometric images of those complex trig identities are, and what their applications are. In this way, you may discover something new and surprising. 

\subsubsection{Double-angle formulas}
Double-angle formulas, by name,  compute the trigonometric functions of $2 \theta$ when the functions of $\theta$ are given. These formulas can be derived from the addition formulas. For example,
\be
\sin(2\theta) = \sin(\theta + \theta) = \sin\theta \cos\theta +  \cos\theta \sin\theta = 2 \sin\theta \cos\theta.
\ee

Similarly, you can derive that 
\be
\cos(2\theta) = \cos(\theta + \theta) = \cos^2\theta -  \sin^2\theta .
\ee

The double-angle formulas for sine and cosine functions are profound. They imply that a high-frequency wave can be represented by waves of half the frequency, since sines and cosines are often used to describe waves and oscillatory motion, where $\theta = 2\pi \omega t$, with $\omega$ the frequency and $t$ time. This powerful formula can lead to very deep mathematical thinking, such as the Fourier series, which states that every smooth curve can be represented as a sum of sines and cosines of different frequencies. 

The importance of the double-angle formula may prompt your curiosity for the triple-angle formula
\bea
\sin(3\theta) &=& \sin(2\theta + \theta) \nonumber \\
 &=& \sin(2\theta) \cos\theta +  \cos(2\theta) \sin\theta 
 \eea
 You can go further to represent this using only $\sin\theta $ and $\cos\theta$. 

\subsection{Trigonometric functions beyond the unit circle}

\begin{figure}[h]
\centering
\includegraphics[width=3.5in]{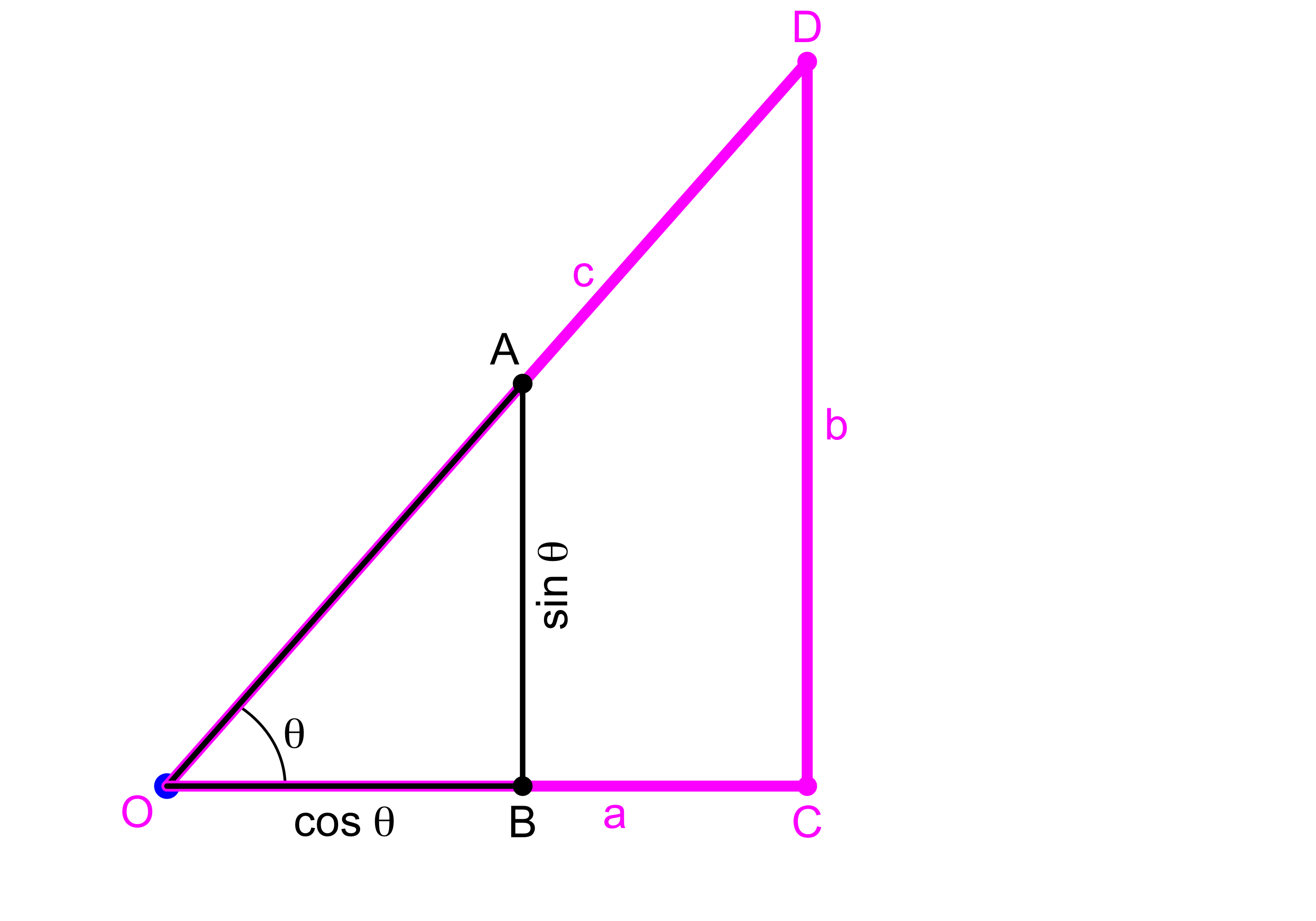}
\caption{Right triangles with different lengths of hypotenuse. }
\label{fig06}
\end{figure}

Right triangles are not limited to the ones with a unit hypotenuse. The hypotenuse can be extended or contracted to any length, as shown in Fig. \ref{fig06}.  

The hypotenuse of the black right triangle $\Delta AOB$ is 1, but that of the magenta right triangle $\Delta COD$ is $c$, which is greater than 1 in the figure, but can be less than 1. These two triangles are similar because their three interior angles are equal. The equal ratios for similar triangles imply that 
\be
\frac{\sin \theta}{1} = \frac{b}{c},
\ee
where $b = \overline{CD}$ is the length of the opposite side of the magenta right triangle $\Delta COD$, and $c = OD$ is the length of its hypotenuse. This equation is commonly used as the definition of $\sin \theta$ in many textbooks. 
Teachers guide students to memorize that sine equals the length of the opposite side divided by the length of the hypotenuse. 

Another set of equal ratios implies 
\be
\frac{\cos \theta}{1} = \frac{a}{c},
\ee
i.e., cosine equals the adjacent side divided by the hypotenuse. 

For math lovers and students with good memory and algebra skills, these definitions have certain advantages, 
particularly from the point of view of analytic geometry, because these are analytic definitions that use formulas to study geometry. Almost all the mathematics teachers are either math lovers or have good algebra skills. These ratios are natural and can be conveniently used to derive many other formulas and effectively solve problems in books. However, as mentioned earlier, most students learn by images in their minds. The algebraic approach is considered abstract and sometimes makes no sense, let alone remembering it and applying it in their careers and lives.

\section{Increments for the sine function}

\begin{figure}[h]
\centering
\includegraphics[width=4.4in]{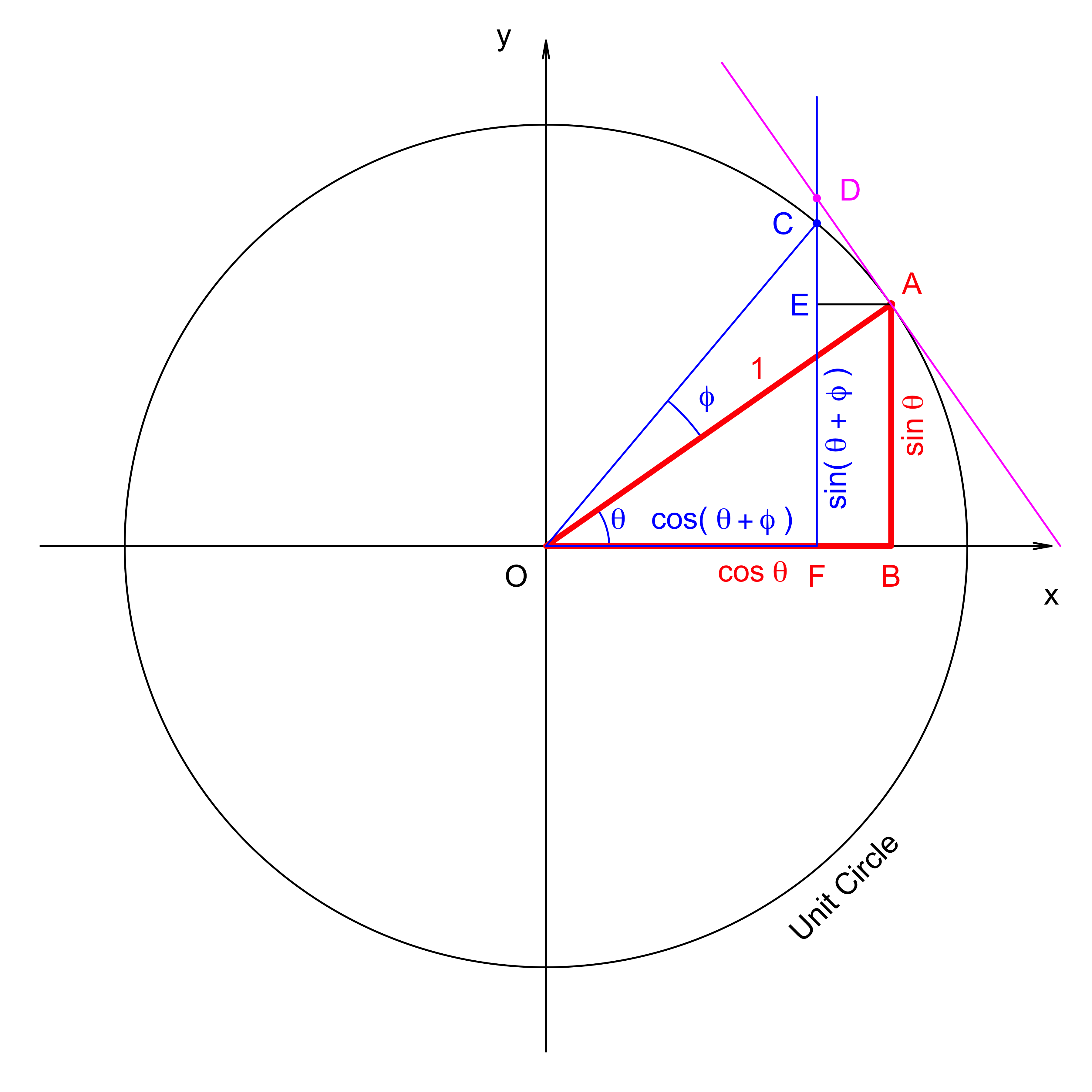}
\caption{Use a figure to derive the derivative formula for the sine function. }
\label{fig07}
\end{figure}

 Most calculus students do not remember or understand the derivation of the sine function's derivative because deriving the derivative formula is an arduous, multi-step process. The process requires several properties of sines and cosines, as shown in most calculus textbooks. Here, we present an intuitive geometric structure for the ratio of increments associated with the derivative.
Hope that our figure (Fig. \ref{fig07}) can leave a lasting impression on students and empower them to explore further by drawing similar or more creative diagrams. We are not providing a rigorous proof of $d(\sin\theta)/d\theta = \cos\theta$,
but we wish to use the figure to show that this derivative formula is geometrically intuitive and reasonable.

Figure \ref{fig07} shows the angle $\theta$, its increment $\phi$ (often denoted by $\Delta \theta$ in textbooks), and their relevant lines and triangles.  Point $A$ corresponds to $\theta$ and $\sin\theta$ on the unit circle, point $C$ to $\theta + \phi$ and $\sin(\theta + \phi)$. The magenta line $AD$ is the tangent line of the circle at point $A$, and is perpendicular to the radial line $OA$. Thus, triangles $\Delta DEA$ and $\Delta OBA$ are similar and provide the following equal ratios:
\be
\frac{DE}{DA} = \frac{OB}{OA}.
\label{eq45}
\ee

In calculus practice, we use a tangent line segment to approximate a segment of a curve or a secant line segment. The small error between the function values of the approximation is one order smaller than 
the small increment of the independent variable. 
Corresponding to Fig. \ref{fig07}, this statement means that the approximation error $CD$
 is in the order of $(\Delta \theta)^2 = \phi^2 = (\stackrel{\frown}{AC}/1)^2$ or smaller, 
 while the increment of the independent variable is $\Delta \theta = \phi = \stackrel{\frown}{AC}/1$, since 
 the length of an arc equals the corresponding interior angle times radius.  Thus, we have the following:
\bea
\overline{OA} &=& 1, \\
\overline{OB} &=& \cos \theta,\\
\overline{DA}  &\approx & \stackrel{\frown}{CA}  = 1\times \phi~~\text{(The Arc Length of} ~AC), \label{eq62}\\
\overline{DE} &\approx & \overline{CE} = \sin(\theta + \phi) - \sin\theta. \label{eq63}
\eea
These equations and Eq. (\ref{eq45}) imply that 
\be
\frac{\sin(\theta + \phi) - \sin\theta}{\phi} \approx \cos \theta .
\label{eq50}
\ee
This suggests that 
\be
\frac{d}{d\theta} \sin \theta = \cos \theta.
\ee

The approximations omit (i) the error $\overline{CD}$ and (ii) the difference of the 
length of the arc $\stackrel{\frown}{AC}$ minus the length of the tangent line segment $\overline{DA}$  when $\phi$ goes to zero. This statement needs a proof, which can be found in most calculus textbooks and is not given here to distract, since our purpose is not to provide a rigorous proof, but to support geometric reasoning and image development in students' minds, with the objective of effective learning. 



\section{Review, explore, and tell }

\subsection{Review our suggested methods of teaching trigonometric functions}
We have learned that the sine function is the length of a half-chord of a unit circle. The chord may be regarded as the string of an archery set when the bow is treated as part of a unit circle (Fig. \ref{fig01}). This concept was developed in India around 500 CE, but the language translation of the term ``half-chord" got corrupted during the spreading process of the mathematical theory from India to the Middle East and then to the Western world. Although the modern function symbol $\sin$, derived from the Latin word ``sinus," is mathematically defined as the Indian half-chord, its literal meaning is not. This may have made it difficult to learn sine functions and may be the root of confusion among students. Some languages, such as Chinese and Japanese, have retained, or partially retained,  the original Sanskrit meaning of "half-chord," despite using the same function symbol $\sin$ as English and other Western languages. 

You may review Figs. \ref{fig01}-\ref{fig03}, hand-draw your own triangles, and mark the definitions of the six trigonometric functions on a figure. You may draw these figures together with your family members and friends and tell them the story of sine. Hope they are amazed, particularly your parents and grandparents. Most likely, they did not learn trigonometry this way. You may try teaching trigonometry to your ten-year-old brother or sister in our way and see if they can get something out of it. 

Note that $\sin \theta$ can determine all the right triangles in \ref{fig03}, but it cannot determine an oblique triangle.  To determine an oblique triangle, another quantity is needed. The corresponding relationship is the law of cosines (see Eq. (\ref{eq24})):
\be
c^2 = 1^2 + a^2 - 2 a\cos \theta.
\ee
Here, $a$ or $c$ is the additional quantity needed. 

Why is the additional quantity needed? This is because the $90^\circ$ angle constraint in a right triangle is removed for an oblique triangle. This removed condition is replaced by $a$ or $c$. 

You can also review the six trigonometric functions and their identities using our provided Python and R code. You may run the code to reproduce the figures in this paper and even generate some variations of these figures. You may use the computer code to verify the identities you can find on the Internet or in books. 

The reviews above may help you develop more images in your mind. These reviews are not based solely on memories but also on your common sense, particularly your visual capabilities. Drawing the figures can also train the synchronization between your mind and hands. Telling the story or teaching others can catalyze your own thinking and motivate you to explore further. 

\subsection{Explore further}

Given $\sin \theta$, every side of all the triangles in Fig. \ref{fig01} can be determined.
Given $\theta$, every angle in Fig. \ref{fig01} can also be determined. In fact, when one side is given for any of the triangles in Fig. \ref{fig01}, all the other sides can be determined. 
Aryabhatta used $\overline{AB} = \sin \theta$ as his primary variable. His choice was the most natural for his astronomical observations, such as observing the Sun relative to the horizon [1]. Can you choose an unnamed side, such as $\overline{BC}$ or $\overline{DE}$, as your primary variable for your observations or measurements? Do you want to explore different options to possibly develop new instruments or machines?

It is commonly agreed that Aryabhata calculated the first table of sine values. He divided the first quadrant into 24 intervals, each of which is 3.75$^\circ$. Then, he used the half-chord definition and developed an algorithm to compute 
$\sin (i \alpha\times \pi/180)$ values  for $i = 1, 2, \cdots, 24$ and $\alpha = 3.75$.
Aryabhata's approach was primarily for calculating planetary positions, as he was better known as an astronomer than as a mathematician [1]. 

Would you like to explore Aryabhatta's original writing [1, 2] and reproduce his sine table following his approach?
Would you like to measure the half-chord length with a ruler or tape? 
Can you use Python or R code to verify your computed or measured values? Mathematical historian and educator, Victor J. Katz, encouraged questions and the use of history in mathematics teaching. In the short article [3], he asked many inspiring questions.
One of them is ``Can we help our students to discover [the Pythagorean theorem]?" Indeed, we can. In calculus teaching, I often guide students to ``discover" different physical laws, such as Newton's second law of motion by data plotting, Faraday's law for 
generating electrical power with a simple magnet-and-wire demo, and Ampere's law when making an electrical motor. 

You may also explore the numerical values of the approximation errors in the process of going from ratios of increments to the derivative of a sine function. Figure \ref{fig07} suggests that the approximation formulas (\ref{eq62}) and (\ref{eq63}) are reasonable and intuitive. 
Would you like to explore these approximations and see how the approximation errors decrease when the increment angle $\phi$ is getting smaller? Namely, for a given $\theta$, you can use Python, R, or a calculator to compute the numerical values of 
$\frac{\sin(\theta + \phi) - \sin\theta}{\phi} $ and $\cos \theta$, and check the differences between the two sides of the approximation formula (\ref{eq50}). Then, you gradually reduce the $\phi$ values and observe the trend of the differences. Do you have confidence to conclude that the limit of the left-hand side of (\ref{eq50}) is equal to $\cos \theta$, as $\phi$ goes to zero?

Of course, infinitely many cases exist in the real or imaginary world relevant to trigonometry for you to explore. The above lists only a few to spark your further thinking. You might invent something of great importance, as Aryabhatta did! 

\subsection{Tell your stories about trigonometry} 

\begin{figure}[h]
\centering
\includegraphics[width=3.0in]{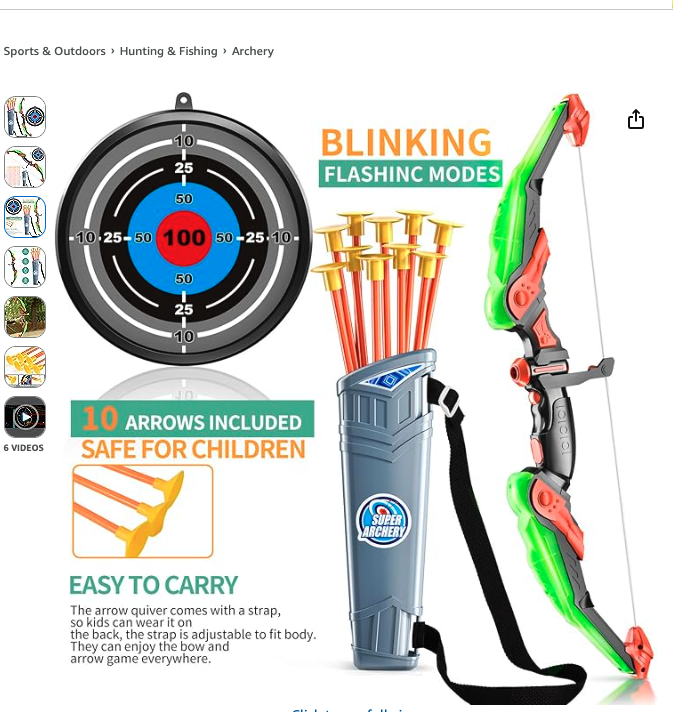}
\includegraphics[width=3.0in]{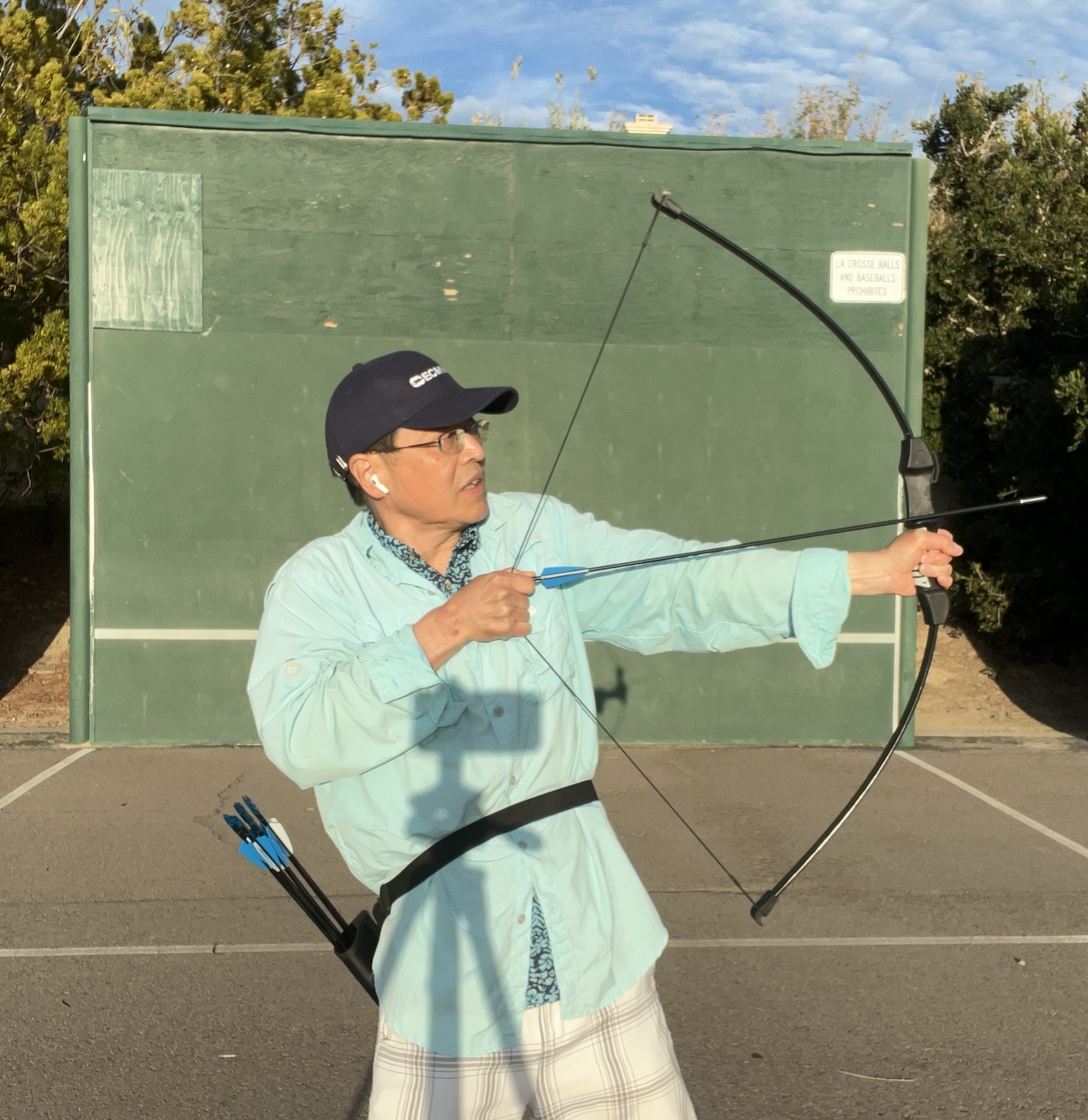}
\caption{The author bought an archery set from AMAZON and tried to feel Aryabhata's archery experience of 17 centuries ago.} 
\label{fig08}
\end{figure}

Would you like to retell the story of India's invention of {\tt jya} and the meaning of sine for the function $\sin x$. 
Would you like to play with an archery set and observe its bow and chord? 
I actually bought a simple toy archery set from AMAZON and played with it many times. See Fig. \ref{fig08}. Of course, 
be safe. A toy archery arrow may still hurt when being hit. So stay away from people when playing with it. 

With friends and at parties, I shared my own sine stories and the experience with my archery toy. Would you like to do the same?

To this point, what images have you got in your mind? What have you experienced? 
Would you like to share your stories on trigonometry? The world wants to hear from you!


\section{Summary and discussion}

We have used the one-figure-facilitates-every-relationship (OFFER) method and the story on the etymology of the mathematics notation $\sin$ to show that all six trigonometric functions can be shown on a single figure. The original name for each function makes sense: $\sin \theta$ is the half-chord, $\tan\theta$ is the length of a tangent line segment ending on the 
$x$-axis, denoted by $C$, and $\sec \theta$ is the half length of a horizontal line segment $\overline{HC}$ that cuts the unit circle into two equal pieces. 
Each complementary function has its own corresponding way: $\cos \theta$ as the other leg of a right triangle with $\sin \theta$, $\cot\theta$ as the extension of the $\tan\theta$ line segment toward the $y$-axis, and $\csc \theta$ as the
length of the line segment on the $y$-axis similar to $\sec \theta$ on the $x-$axis. This OFFER method for describing fundamental mathematical concepts may help reduce misunderstandings caused by the conventional algebra-based definitions of trigonometric functions.  As discussed earlier, our OFFER approach may help people avoid being puzzled by the following two expressions,
\be
\sec \theta = \frac{1}{\cos \theta}, ~~~ \csc \theta = \frac{1}{\sin \theta} 
\ee
when these expressions are treated as consequences instead of definitions. 


Our paper is written based on the pedagogy of story-picture-observe-review-tell (SPORT), in order to make learning effective.  This mathematics education pedagogy is a concrete practice of learning-by-doing. Another Chinese saying, also attributed to Confucius, is that ``I hear, I forget; I see, I remember; and I do, I understand."  Many mathematics instructors and university professors are frustrated that students forget what they learned even a year ago. Why? 
Professors often attribute this to students' negligence and indolence. This attribution might be unfair! The names of these concepts are not encountered in daily life, so students cannot form images in their minds.  Thus, it is not because the concepts are hard, but the names do not make sense. Examples are numerous, such as functions, algebra, calculus, and eigenvectors.  Those mathematical concepts are actually very easy and intuitive; however, they are taught in a very hard way, which consequently makes students not even want to remember them.  

In Japanese, function is \begin{CJK}{UTF8}{min}関数\end{CJK}, pronounced as [{\tt kansuu}] and literally meaning related numbers or enclosed numbers.  In Chinese and Japanese, algebra is \begin{CJK*}{UTF8}{gbsn}(代数)\end{CJK*}, meaning substituting numbers. In German, calculus is {\it infinitesimalrechnung}, meaning infinitesimal accounting, and eigenvector is {\it eigenvektor}, meaning my own vector, a special vector whose orientation is not changed by a matrix multiplication [4, 5]. These mathematical names in the said languages make good sense, but may not in English. For instance, many American students do not know the meaning of the German or Dutch word ``eigen", meaning one's own or self. The term ``eigenvector" was coined by German mathematician David Hilbert as ``eigenvektor" at the beginning of the 1900s and was brought to the American scientific communities by his student Richard Courant and John von Neumann soon after that [6]. 

The Chinese name of ``Sine function"  as {\tt zhengxian hanshu}  \begin{CJK*}{UTF8}{gbsn}(正弦函数)\end{CJK*} makes good sense to Chinese students, since they can form a picture in their mind of a chord directly opposite to an angle. In contrast, the English name sine is completely disjoint from the geometric meaning of $\sin(x)$. It is very hard for American students to form a picture of sine from the name alone. They are often forced to remember the formula definition of sine. This ``forced" process is a typical learning approach from memory, in contrast to a more effective method of learning from common sense. 

Although Chinese students have the advantage of learning from common sense when discussing the sine function, the notation $\sin(x)$ can still confuse them, as they also pronounce $\sin(x)$ the same as English, [{\tt sine x}]. This pronunciation makes no sense to many of them. 

Would it be helpful if our instructors and textbooks could explain intrinsic mathematical concepts in words or stories relevant to students' daily lives? The etymology of some words, like sine, as shown in this paper, can be a good story and help students build mental pictures. 

For the benefit of not only English-speaking students but also others, can we name
$\text{hch}(x)$ for $\sin(x)$, and $\text{coc}(x)$ for $\cos(x)$? Pronounce $\text{hch}(x)$ as half-chord $x$ in English, and [{\tt banxuan x}] in Chinese. 
Pronounce $\text{coc}(x)$ as co-chord $x$ in English, and [{\tt yuxuan x}] in Chinese. The latter has been used in Chinese pronunciation for $\cos(x)$ anyway, although most Chinese students often pronounce it the English way, [{\tt cosine x}].  

We have no intention to suggest the change of hundreds of years of traditional notation of $\sin(x)$ and $\cos(x)$. Our discussion here is to inspire students to think more deeply and explore further.  

The mathematics education pedagogy and learning strategy suggested in this paper may be summarized as the following four actions:
\begin{description}
\item{(a)} Do SPORT to practice effective learning, 
\item{(b)} Make an OFFER to learn a series of concepts using one figure,
\item{(c)}  Learn RUM materials first, and  
\item{(d)} Leave OUT materials for future studies. 
\end{description}

Our final advice: {\it SPORT OFFERs RUM, and don't go OUT!}

\section{Code availability}
The R and Python code created by the author with AI assistance for this paper are in a .zip file and freely available at GitHub for download\\
{\tt https://github.com/SamuelSPShen/OFFER-for-Trig-Functions}\\
Also see Ref. [7]:  {\tt https://doi.org/10.5281/zenodo.21267201}

\vskip 0.2cm

\noindent  {\bf Acknowledgments} ~~ This study was supported in part by the US National Science Foundation [\#IIS2324008] and 
U.S. National Oceanographic and Atmospheric Administration [\#NA22SEC4810016].

\vskip 0.2cm
\noindent  {\bf Disclosure statement } ~~ No potential conflict of interest was reported by the author.

\vskip 0.2cm
\noindent  {\bf Declaration of generative AI use} ~~ AI was used to (i) generate Fig. \ref{fig02}, (ii) convert part of the R code into Python for this paper, and (iii) make some online searches for wording, grammatical expressions, and historical materials. 

\vskip 0.5cm

\noindent {\Large REFERENCES}

\begin{description}

\item{[1]}  Singh V A, Kumar A.  Aryabhata and the Construction of the First Trigonometric Table. {\it arXiv:2309.13577}, 2023.

\item{[2]}  Shukla K S, K V Sarma.  Aryabhatiya of Aryabhata. Indian National Science Academy, 1976.

\item{[3]} Katz V J. Some ideas on the use of history in the teaching of mathematics. {\it For the Learning of Mathematics}, 17, 62-63, 1997.

\item{[4]}.	Shen S S P, Somerville R C J. {\it Climate Mathematics: Theory and Applications}, Cambridge University Press, 2019.

\item{[5]}.	Shen S S P,  North G R.  {\it Statistics and Data Visualization in Climate Science with R and Python}, Cambridge University Press, 2023.

\item{[6]} Tou E R. Math origins: Eigenvectors and eigenvalues. Convergence. 15, 12, 2018. DOI: 10.4169/convergence20181126. 

\item{[7]} Shen, S. S.P .  Python and R Code for "OFFER (One-Figure-Faciliates-Every-Relationship) for Trigonometric Functions" by Samuel S.P. Shen (V1.0.0). Zenodo, , 2026. \\{\tt https://doi.org/10.5281/zenodo.21267201}
\end{description}

\vskip 0.5cm

\noindent {{\bf  Summary}} This article describes an approach of using a single figure to explain all six trigonometric functions, their fundamental identities, and the cosine law. We name this method OFFER (one-figure-facilitates-every-relationship) for the convenience of mathematics teaching. The trigonometric identities and formulas are shown to have geometric interpretations using the corresponding triangles, whether right or oblique.  The article also describes the SPORT (pedagogy of the story-picture-observe-review-tell) approach to effective mathematics learning. Stories are told about the invention of the sine function by the Indian mathematician Aryabhatta, the English term``sine", and an archery toy. It is suggested that students first learn the RUM (relevant, useful, and modern) materials and review them before spending too much time on the OUT (only-useful-for-teaching) materials. 

\vskip 0.5cm

\noindent {\bf SAMUEL S.P. SHEN} Distinguished Professor of Mathematics and Statistics at San Diego State University, and Visiting Research Mathematician at Scripps Institution of Oceanography, University of California – San Diego. Formerly, he was McCalla Professor of Mathematical and Statistical Sciences at the University of Alberta, Canada, and President of the Canadian Applied and Industrial Mathematics Society. He has held visiting positions at the NASA Goddard Space Flight Center, the NOAA Climate Prediction Center, and the University of Tokyo. Shen holds a B.Sc. degree in Engineering Mechanics and a Ph.D. degree in Applied Mathematics.

\end{document}